\documentclass[12pt]{amsart}
\usepackage{amssymb,amsmath}
\usepackage{graphicx}
\usepackage{colordvi}

\numberwithin{equation}{section}

\newtheorem{theorem}[equation]{Theorem}

\newtheorem{proposition}[equation]{Proposition}
\newtheorem{corollary}[equation]{Corollary}

\theoremstyle{definition}
\newtheorem{definition}[equation]{Definition}

\newtheorem{example}[equation]{Example}

\newtheorem{remark}[equation]{Remark}
\newcounter{FNC}[page]
\def\fauxfootnote#1{{\addtocounter{FNC}{2}$^\fnsymbol{FNC}$%
     \let\thefootnote\relax\footnotetext{$^\fnsymbol{FNC}$#1}}}
\newcommand{\rat}{\ {\relbar\rightarrow}\ }
\newcommand{\longrat}{{\;\relbar\relbar\rightarrow\;}}

\newcommand{\zwei}[2]{\left[\!\begin{array}{r}\!#1\!\\\!#2\!\end{array}\!\right]}
\newcommand{\kzwei}[2]{\left[\substack{#1\\#2}\right]}
\newcommand{\smi}{\hspace{0.6em}}

\newcommand{\ba}{{\bf a}}
\newcommand{\bb}{{\bf b}}
\newcommand{\bc}{{\bf c}}
\newcommand{\bd}{{\bf d}}
\newcommand{\be}{{\bf e}}
\newcommand{\bp}{{\bf p}}

\newcommand{\bx}{{\bf x}}
\newcommand{\bv}{{\bf v}}

\newcommand{\calA}{{\mathcal A}}
\newcommand{\calB}{{\mathcal B}}

\newcommand{\C}{{\mathbb C}}
\renewcommand{\P}{{\mathbb P}}
\newcommand{\R}{{\mathbb R}}
\newcommand{\Z}{{\mathbb Z}}

\newcommand{\QED}{\hfill\raisebox{-5pt}{\includegraphics[height=14pt]{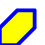}}\vspace{6pt}}

\begin{document}

\title{Linear precision for parametric patches}

\author{Luis David Garcia-Puente}
\address{Department of Mathematics and Statistics\\
         Sam Houston State University\\
         Huntsville\\
         TX \ 77341\\
         USA}
\email{lgarcia@shsu.edu}
\urladdr{http://www.shsu.edu/\~{}ldg005}
\author{Frank Sottile}
\address{Department of Mathematics\\
         Texas A\&M University\\
         College Station\\
         TX \ 77843\\
         USA}
\email{sottile@math.tamu.edu}
\urladdr{http://www.math.tamu.edu/\~{}sottile}
\thanks{Work of Sottile supported by NSF CAREER grant DMS-0538734,
        by the Institute for Mathematics and its Applications with funds provided by the
	National Science Foundation, and by Peter Gritzmann of the Technische
         Universit\"at M\"unchen.} 
\begin{abstract}
%
 We give a precise mathematical formulation for the notions of a parametric patch and linear
 precision, and establish their elementary properties.
 We relate linear precision to the geometry of a particular linear projection, giving
 necessary (and quite restrictive) conditions for a patch to possess linear precision.
 A main focus is on linear precision for Krasauskas' toric patches, which we show is
 equivalent to a certain rational map on $\C\P^d$ being a birational isomorphism.
 Lastly, we establish the connection between linear precision for toric surface patches
 and maximum likelihood degree for discrete exponential families in algebraic statistics,
 and show how iterative proportional fitting may be used to compute toric patches.
\end{abstract}

%
\maketitle

\section*{Introduction}
B\'ezier curves and surfaces are the fundamental units in geometric modeling.
There are two basic shapes for surfaces---triangular B\'ezier patches
and rectangular tensor product patches. 
Multi-sided patches are needed for some applications, and there are 
several control point schemes for $C^\infty$ multi-sided patches.
These include the $S$-patches of Loop and DeRose~\cite{LdR89}, 
Warren's hexagon~\cite{Wa92}, Kar{\v{c}}iauskas's $M$-patches~\cite{Ka03},
and the multi-sided toric B\'ezier patches of
Krasauskas~\cite{KR02}. 
(Relationships between these and other patches are discussed
in~\cite{KK00}.)
Parametric patches are general control point schemes for $C^\infty$
patches whose shape is a polygon or polytope.
They include the patch schemes just mentioned, as well as 
barycentric coordinates for polygons and polytopes~\cite{Fl,Wa75,Wa96}.

The success and widespread adoption of B\'ezier and tensor-product patches is
due in part to their possessing many useful mathematical properties.
Some, such as affine invariance and the convex hull property, are
built into their definitions and also hold for the more general
parametric patches.
Other properties, such as de Castlejau's algorithm 
for computing B\'ezier patches, come from the specific form of their 
Bernstein polynomial blending functions.
Linear precision is the ability of a parametric patch to replicate linear
functions. 
When the blending functions of a parametric patch correspond to the
vertices of a polytope, these blending functions give barycentric coordinates precisely
when the patch has linear precision. 
In this way, blending functions for a parametric patch having linear precision
are barycentric coordinates for general control point schemes. 

For us, a 
(parametric) \Blue{{\it patch}} is a collection of non-negative blending 
functions indexed by a finite set $\calA$ of points in $\R^d$, 
where the common domain of the blending functions is the convex hull $\Delta$
of $\calA$. 
A collection of control points in $\R^\ell$ indexed by $\calA$ is used to
define a map from $\Delta$ to $\R^\ell$.
The blending functions determine the internal structure of this map
and the basic shape, $\Delta$, of its image, while the control points 
determine how the image lies in $\R^\ell$.
Choosing the control points to be the points of $\calA$ gives the 
tautological map, and the patch has linear precision when this 
tautological map is the identity on $\Delta$.
We show that every patch has a unique reparametrization (composing the
blending functions with a homeomorphism of $\Delta$) having linear
precision (Theorem~\ref{T:Lin_Prec}).
This generalizes Theorem~8.5 of~\cite{So03}, which
was for toric patches.

The blending functions of a patch may be arbitrary
non-negative $C^\infty$ functions. 
A patch is rational if it has a reparametrization having polynomial 
blending functions.
For a rational patch, a map given by control points corresponds to a
linear projection of a projective algebraic variety associated to its
blending functions. 
Its unique reparametrization having linear precision has rational blending
functions if and only if a certain canonical map defined on
this variety is a birational isomorphism (Theorem~\ref{Th:LP=birational}),
which implies that this variety has a maximally degenerate position with
respect to a canonical linear subspace given by the set $\calA$ (Theorem~\ref{Th:Necessary}).

We apply this analysis to Krasauskas's toric patches, which are rational.
The shape of a toric patch is a pair $(\calA,w)$ where $\calA$ is a set
of integer points in $\Z^d$ and $w$ is a collection of positive real numbers indexed by
$\calA$. 
We show that iterative proportional fitting~\cite{DR72}, a simple numerical algorithm from
statistics, computes the blending functions that  have linear precision. 
This algorithm was suggested to us by Bernd Sturmfels.
It may form a basis for algorithms to manipulate these patches.

Krasauskas~\cite[Problem 3]{Kr06} asked whether any toric patches (besides the classical B\'ezier
simploids~\cite{dRGHM}) admit a rational reparametrization having linear precision. 
A toric patch of shape $(\calA,w)$ corresponds to a homogeneous polynomial $f=f_{\calA,w}$ 
whose dehomogenization is the sum of monomials with exponents $\calA$ and coefficients $w$.
The toric patch admits a rational reparametrization having linear precision if and only if the toric
differential 
\[
   D_{\mbox{\scriptsize\rm toric}}f\ :=\ 
   \Big[x_0\tfrac{\partial}{\partial x_0}f\;\colon\;
    x_1\tfrac{\partial}{\partial x_1}f\;\colon\;\dotsb
    \;\colon\;x_d\tfrac{\partial}{\partial x_d}f\Big]\ 
    \colon\ \C\P^d\ \longrat\ \C\P^d
\]
defines a birational isomorphism (Corollary~\ref{Cor:Reformulation}).
This analysis of linear precision is used in~\cite{BRS}
to classify which toric surface patches can have linear precision.

In Krasauskas's question about linear precision, he allowed the points $\calA$ indexing
the blending functions of toric patches to move within their (fixed) convex hull, 
keeping the same blending functions (See Example~\ref{Ex:Pentagon}).
Our analysis of linear precision for general parametric
patches will help to address that version of his question.
\bigskip

In Section 1, we define parametric patches  and show
that every patch has a unique reparametrization that has linear
precision.
In Section 2 we show that if a
rational patch has linear precision, then an algebraic variety we
obtain from the blending functions has an exceptional 
position with respect to a certain linear space.
We define toric patches in Section 3, and in Section 4 we 
explain how iterative proportional fitting computes the blending functions for toric
patches that have linear precision.

%
%
%
\section{Linear precision for parametric patches}\label{S:PP}

We interpret the standard definition of a mapping via control
points and blending functions (see for example \cite[\S 2]{KK00}) 
in a general form convenient for our discussion.
All functions here are smooth ($C^\infty$) where defined and
real-valued unless otherwise stated.
Let $\R_>$ be the set of strictly positive real numbers and $\R_\geq$ the
set of  non-negative real numbers.

Let $\calA$ be a finite set of points in $\R^d$, which we shall use as
geometrically meaningful indices.
A control point scheme for parametric patches, or 
\Blue{({\it parametric}) {\it patch}} for short, is a collection
$\beta=\{\beta_\ba\mid\ba\in\calA\}$ of non-negative functions, called
\Blue{{\it blending functions}}.
The common domain of the blending functions is the convex hull $\Delta$ of
$\calA$, which we call the \Blue{{\it domain polytope}}.
We will always assume that $\Delta$ is full dimensional in that it has dimension $d$.
We assume that the blending functions do not vanish simultaneously 
at any point of $\Delta$.
That is, the blending functions have no base points in their domain.

A set $\{\bb_\ba\mid\ba\in\calA\}\subset\R^\ell$ of \Blue{{\it control points}}
gives a map  $F \colon\Delta\to\R^\ell$ defined by
 \begin{equation}\label{Eq:P_Patch}
   F(x)\ :=\ 
   \frac{\sum_{\ba\in\calA} \beta_\ba(x)\,\bb_\ba }
        {\sum_{\ba\in\calA} \beta_\ba(x)}\ .
 \end{equation}
The denominator in~\eqref{Eq:P_Patch} is positive on $\Delta$ and so
the map $F$ is well-defined. 

\begin{remark}\label{R:Comments}
  Positive weights, or shape parameters, 
  $\{w_\ba\in\R_>\mid\ba\in\calA\}$ scaling the blending functions are often included
  in the definition of the map $F$.
  We have instead absorbed them into our notion of blending functions.
  See Remarks~\ref{R:Restrictive} and \ref{Rem:weights}.\hfill\QED
\end{remark}

The blending functions $\{\beta_\ba\mid\ba\in\calA\}$ are 
\Blue{{\it normalized}} if they form a partition of unity,
\[
  \sum_{\ba\in\calA} \beta_\ba(x)\ =\ 1\,.
\]
As the denominator in~\eqref{Eq:P_Patch} is (strictly) positive on $\Delta$, we may divide
each blending function by this denominator to obtain normalized blending functions.
\[
  \mbox{Redefine}\quad
  \beta_\ba\ :=\ \frac{\beta_\ba}{\sum_{\ba\in\calA} \beta_\ba}\,.
   \qquad\mbox{Then}\quad
   \beta_\ba(x)\ \geq\ 0
   \quad\mbox{and}\quad
  \sum_{\ba\in\calA} \beta_\ba(x)\ =\ 1\,,
\]
for $x\in\Delta$.
For normalized blending functions, the formula~\eqref{Eq:P_Patch} becomes
 \begin{equation}\label{Eq:Conv_Combin}
   F(x)\ =\ \sum_{\ba\in\calA} \beta_\ba(x) \bb_\ba\,.
 \end{equation}
We deduce two fundamental properties of parametric patches.
\begin{itemize}
 \item {\bf Convex hull property.} 
    For $x\in\Delta$, $F(x)$ is a convex combination of the control
    points.
    Thus $F(\Delta)$ lies in the convex hull of the control points.
 \item {\bf Affine invariance.} If $\Lambda$ is an affine function on
   $\R^\ell$, then  $\Lambda(F(x))$ is the map with the same blending
   functions, but with control points
   $\{\Lambda(\bb_\ba)\mid\ba\in\calA\}$.  
\end{itemize}

A patch is \Blue{{\it non-degenerate}} if its blending functions are
linearly independent.
This implies in particular that $F(\Delta)$ cannot collapse into a point
unless all the control points are equal.

 We make our key definition.

\begin{definition}\label{DE:LP}
 A patch 
 $\{\beta_\ba\mid\ba\in\calA\}$ has \Blue{{\it linear precision}} if 
 for every affine function $\Lambda(x)$ defined on $\R^d$, 
\[
   \Lambda(x)\ =\    
     \sum_{\ba\in\calA} \beta_\ba(x)\,\Lambda(\ba)\,.
\]
 That is, if the patch can replicate affine functions.
 This notion depends strongly on the positions of the points in
 $\calA$. \hfill\QED
\end{definition}

The \Blue{{\it tautological map $\tau$}} of a given patch is the
map~\eqref{Eq:Conv_Combin} when the control points are taken
to be the corresponding points of $\calA$,
 \begin{equation}\label{Eq:Tautological}
   \tau(x)\ :=\    
     \sum_{\ba\in\calA} \beta_\ba(x)\,\ba\,.
 \end{equation}
By the convex hull property, $\tau(\Delta)\subset \Delta$. 
By affine invariance, the patch has linear precision if and
only if $\tau$ is the identity function, $\tau(x)=x$.
We record this fact.

\begin{proposition}\label{P:LP_criterion}
  A patch has linear precision if and only if its tautological
  map is the identity map on the domain polytope $\Delta$.
\end{proposition}

\begin{remark}\label{R:Restrictive}
 This is a more restrictive notion of linear precision
 than is typically considered.
 It is more common to define a patch to be a collection 
 $\{\beta_i\mid i=1,\dotsc,n\}$ of
 blending functions, and this has linear precision if we have
\[
   x\ =\ \frac{\sum_{i=1}^n w_i\,\beta_i(x)\, \ba_i}
    {\sum_{i=1}^n w_i\,\beta_i(x)}\ ,
\]
 for some non-negative \Blue{{\it weights}} $\{w_i\mid i=1,\dotsc,n\}$ and
 some points $\calA=\{\ba_i\mid i=1,\dotsc,n\}$ whose convex hull is the  
 domain polytope $\Delta$ of the blending functions.

 Our more restrictive definition of patches and linear precision, where we
 incorporate the weights and the points $\calA$ into our definition of blending functions,  
 allows us to derive precise criteria with which to study linear precision.
 Our intention is to employ these criteria to study patches in the
 generality in which other authors had worked. \hfill\QED
\end{remark}

\begin{example}\label{Ex:BC}
 Fix an integer $n>0$ and let $\calA:=\{\frac{i}{n}\mid \,0\leq i\leq n\}$
 be the set of $n{+}1$ equally distributed points in the unit interval.
 The $i$th \Blue{{\it Bernstein polynomial\/}} is $\beta_i:=\binom{n}{i}x^i(1{-}x)^{n-i}$,
 which we associate to the point $\frac{i}{n}\in\calA$.
 These form a partition of unity,
\[
     \sum_{i=0}^n \beta_i(x)\ =\ 
     \sum_{i=0}^n \binom{n}{i}x^i(1{-}x)^{n-i}\ =\ 
      \bigl(x\, +\, (1{-}x)\bigr)^n\ =\ 1\,,
\]
 and are therefore normalized blending functions.
 Given control points $\bb_{i}\in\R^3$, formula~\eqref{Eq:P_Patch} 
 becomes
\[
      F(x)\ =\ 
   \sum_{i=0}^n \binom{n}{i}x^i(1-x)^{n-i} \bb_i \,,
\]
 which is the classical formula for a B\'ezier curve of degree $n$ in $\R^3$.

 This patch has linear precision.
 First, note that $\frac{i}{n}\cdot\binom{n}{i}= \binom{n-1}{i-1}$.
 Then $\tau(x)$ is 
 \begin{eqnarray*}
  \sum_{i=0}^n\binom{n}{i}x^i(1-x)^{n-i}\cdot \frac{i}{n}
  &=&
  \sum_{i=1}^n\binom{n-1}{i-1}x^i(1-x)^{n-i}\\
   &=& x\cdot 
    \sum_{j=0}^{n-1}\binom{n-1}{j}x^j(1-x)^{n-1-j}
  \quad=\quad x\,.\makebox[.1cm][l]{\hspace{36pt}$\QED$}
 \end{eqnarray*}
\end{example}

\begin{example}\label{Ex:Wach}
 Let $\calA$ be the vertices of the hexagon, given
%
%
 below by its defining inequalities.
%
 \[
  \begin{picture}(250,172)(-51,-11)
   \put( 0, 0){\includegraphics[height=150pt]{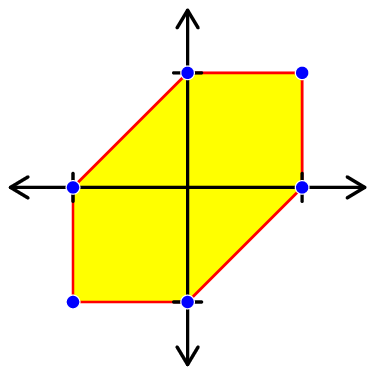}}
   \put(120.5,58){$1$}     \put(60,120){$1$}
   \put(83,23.5){$-1$}   \put(13.8,84){$-1$}

                            \put(100,153){$1-y>0$}
   \put(-35,122){$1+x-y>0$} \put(148, 95){$1-x>0$}
   \put(-51, 47){$1+x>0$}   \put(113, 22){$1+y-x>0$}
   \put(  2,-10){$1+y>0$}

                                     \put(111.5,148){\vector(-1,-2){9.5}}
   \put(31.5,118.5){\vector(1,-1){14.5}}\put(146, 98){\vector(-1, 0){20}}
   \put( 1  ,50){\vector(1, 0){20}}  \put(118.5, 31.5){\vector(-1, 1){14.5}}
   \put(38.5, 2){\vector(1, 2){9.5}}
  \end{picture}
 \]
Let the blending function associated to a vertex be the product of the linear forms
defining edges that  do not contain that vertex.
 \begin{equation}\label{Eq:WHex}
  \begin{array}{rcl}
   \beta_{\kzwei{-1}{-1}}&:=&(1+x-y)(1-y)(1-x)(1+y-x)\\\rule{0pt}{14pt}
   \beta_{\kzwei{-1}{\smi 0}} &:=&(1-y)(1-x)(1+y-x)(1+y) \\\rule{0pt}{14pt}
   \beta_{\kzwei{0}{1}}  &:=&(1-x)(1+y-x)(1+y)(1+x) \\\rule{0pt}{14pt}
   \beta_{\kzwei{1}{1}}  &:=&(1+y-x)(1+y)(1+x)(1+x-y) \\\rule{0pt}{14pt}
   \beta_{\kzwei{1}{0}}  &:=&(1+y)(1+x)(1+x-y)(1-y) \\\rule{0pt}{14pt}
   \beta_{\kzwei{\smi 0}{-1}} &:=&(1+x)(1+x-y)(1-y)(1-x)
 \end{array}
\end{equation}
The normalized blending functions have denominator
$\varphi(x,y):=2(3+xy-x^2-y^2)$, which is strictly positive on  
$\Delta$.
These normalized blending functions have linear precision.
\[
   \sum_{\ba\in\calA} \beta_{\ba}(x,y)\,\ba\ =\ 
    \varphi(x,y)  \zwei{x}{y} \ .
\]
In fact, these normalized blending functions are
Wachspress's barycentric coordinates~\cite{Wa75} for this hexagon.\hfill\QED
\end{example}

When $\calA$ is the set of vertices of  $\Delta$
(so that $\calA$ is in convex position), then normalized blending
functions of a patch with linear precision are 
\Blue{{\it barycentric coordinates}} for $\Delta$.
By this we mean in the sense of~\cite{Fl,Wa75,Wa96}: nonnegative functions indexed by the
extreme points of $\Delta$ that have linear precision in the sense of
Definition~\ref{DE:LP}.\smallskip 

A \Blue{{\it reparametrization}} of a patch 
$\{\beta_\ba\mid\ba\in\calA\}$ by a  homeomorphism 
$\psi\colon \Delta\to\Delta$ 
is a new patch with blending functions
$\beta_\ba\circ \psi$ for $\ba\in\calA$.
A patch is \Blue{{\it proper}} if its tautological
map~\eqref{Eq:Tautological} is a homeomorphism.
This condition is necessary for the patch to have linear precision.
If a patch is proper, then reparametrizing it by 
the inverse $\tau^{-1}$ of $\tau$ gives a new patch whose
tautological map is $\tau\circ\tau^{-1}$, the identity
function, and $\tau^{-1}$ is the unique homeomorphism of $\Delta$ with
this property. 
We record this straightforward, yet fundamental result about linear precision. 

\begin{theorem}\label{T:Lin_Prec}
  A proper patch has a unique reparametrization that has linear precision.
\end{theorem}

Theorem~\ref{T:Lin_Prec} suggests that it will be fruitful to discuss
patches up to reparametrization.
In Section~\ref{S:Geometry}, we show that a natural geometric object
$X_\beta$ associated to a patch is invariant under 
reparametrizations and thus represents the patch up to reparametrization.

A patch is \Blue{{\it rational}} if it has a reparametrization whose
normalized blending functions are rational functions (quotients of
polynomials).
A patch has \Blue{{\it rational linear precision}} if
it is proper and its reparametrization having linear precision has blending
functions that are rational functions.
Such a patch is necessarily rational.
We seek criteria that determine when a rational patch has
rational linear precision.

Theorem~\ref{T:Lin_Prec} and our preceding discussion concerns
reparametrizations of a patch having linear precision.
We may also alter a patch by multiplying its blending functions by
positive real numbers, typically called weights, or  
by moving the points of $\calA$ that are not vertices of the domain
polytope $\Delta$. 
This second change does not alter any map~\eqref{Eq:P_Patch} given by
control points, but it will change the tautological map
$\tau$~\eqref{Eq:Tautological}. 
Our analysis in the next section may be helpful in addressing 
whether it is possible to tune a given patch (using weights or
moving the points in $\calA$) to obtain one that has 
rational linear precision. 
Example~\ref{Ex:Pentagon} (due to Krasauskas and Kar{\v{c}}iauskas)
shows how one toric patch may be tuned to achieve linear precision.

%
%
%
\section{The geometry of linear precision}\label{S:Geometry}

We introduce an algebraic-geometric formulation of patches to clarify the
discussion in Section~\ref{S:PP} and to provide tools with which to
understand and apply Theorem~\ref{T:Lin_Prec}.
This leads us to discuss the relevance for geometric modeling of the
subtle difference between rational varieties and unirational varieties.
Lastly, we use this geometric formulation to give a geometric
characterization of when a patch has rational linear precision.
We recommend the text~\cite{CLO} for additional background on algebraic
geometry. 
We first review linear projections, which are the geometric counterpart
of control points.

\subsection{Linear projections}

We consider $\R^\ell$ as a subset of the 
$\ell$-dimensional real projective space $\R\P^\ell$ via the embedding 
\[
  \R^\ell\ \ni\ z\ \longmapsto [1,z]\ \in\ \R\P^\ell\,.
\]
A point $[z_0,z_1,\dotsc,z_\ell]\in\R\P^\ell$ lies in this copy of
$\R^\ell$ if and only if $z_0\neq 0$.
In that case, the corresponding point of $\R^\ell$ is
 \begin{equation}\label{Eq:CoordsRd}
   \Bigl( \frac{z_1}{z_0},\, \frac{z_2}{z_0}, \,\dotsc,\,
       \frac{z_\ell}{z_0}\Bigr)\,.
 \end{equation}

Let $\calA\subset\R^d$ be a finite set of points and
write $\R\P^\calA$ for the real projective
space whose coordinates are indexed by $\calA$.
Let $b=\{\bb_\ba\in\R^\ell\mid \ba\in\calA\}$ be a collection of control
points, which we regard as  points in $\R\P^\ell$ so that
$b=\{[1,\bb_\ba]\in\R\P^\ell\mid \ba\in\calA\}$.

Given a point $y=[y_\ba\mid\ba\in\calA]\in \R\P^\calA$, if the sum
 \begin{equation}\label{Eq:sum}
   \sum_{\ba\in\calA}  y_\ba \cdot(1, \bb_\ba) \ \in \ \R^{\ell+1}
 \end{equation}
is non-zero then it represents a point
$\pi_b(y)$ in $\R\P^\ell$.
This map $\pi_{b}$ is a \Blue{{\it linear projection}} 
 \begin{equation}\label{Eq:Lin_Proj}
  \pi_b\ \colon\ \R\P^\calA\ \longrat\ \R\P^\ell\,.
 \end{equation}
We use a broken arrow $\rat$ to indicate that $\pi_b$ is defined only
on the complement of the set 
$E_b\subset\R\P^\calA$ where the sum~\eqref{Eq:sum} vanishes.
This linear subspace is the \Blue{{\it center of projection}}.

When the control points affinely span $\R^\ell$,
the center of projection $E_b$ has codimension $\ell\!+\!1$ in $\R\P^\calA$.
The inverse image $\pi_b^{-1}(x)$ of a point $x\in\R\P^\ell$
is the set $H\setminus E_b$, where $H\subset\R\P^\calA$ is a
linear subspace of codimension $\ell$ containing $E_b$.
If $L$ is any $\ell$-dimensional linear subspace of $\R\P^\calA$ that
does not meet $E_b$, then $\pi_b$ maps $L$ isomorphically to 
$\R\P^\ell$.
Identifying $L$ with $\R\P^\ell$ gives an explicit
description of the map $\pi_b$.
If $y\in\R\P^\calA\setminus E_b$, then 
\[
  \pi_b(y)\ =\ L\ \cap\ \left(\overline{ y,\,E_b }\right)\,,
\]
the intersection with $L$ of the linear span of $y$ and the center
$E_b$.  
If $y\in L$, then its inverse image under $\pi_b$ is 
$(\overline{ y,\,E_b })\setminus E_b$.

\begin{figure}[htb]
 \[
  \begin{picture}(210,150)(0,-2)
   \put(0,0){\includegraphics{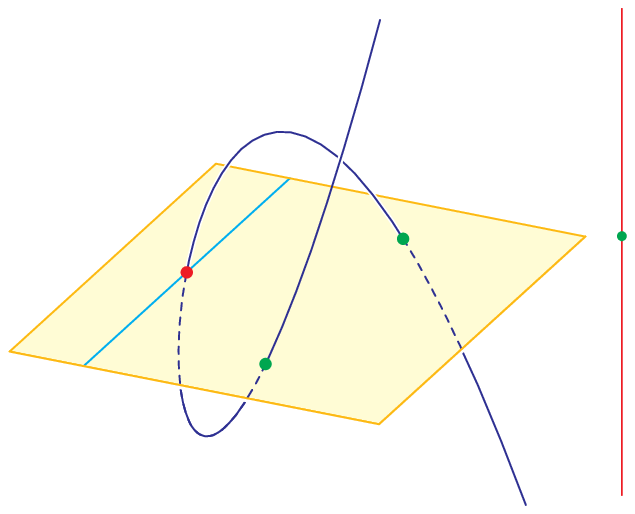}} 

   \put(10,130){$\pi$}\put(20,133){\vector(1,0){20}}

   \put(0,62){$E$}\put(12,64){\vector(2,-1){20}}
   \put(18,79){$B$}\put(30,82){\vector(2,-1){20}}
   
   \put(135,114){$\overline{ y,\, E}$}
   \put(145,109){\vector(-1,-4){5}}

   \put(207,79){$y$}\put(204,81){\vector(-1,0){20}}
   \put(207,119){$L$}\put(204,121){\vector(-1,0){20}}

   \put(155,8){\Blue{$C$}}

   \put(100,12.5){\vector(-3, 4){21}}
   \put(116.5,12.5){\vector(0,1){62}}
   \put(92,-1){$\pi^{-1}(y)$}

  \end{picture}
 \]
\caption{A linear projection $\pi$ with center $E$.\label{F:cubic_Projection}}
 
\end{figure}

For example, Figure~\ref{F:cubic_Projection}
shows the effect of a linear projection $\pi$
on a cubic space curve $C$.
The center of projection is a line, $E$, which meets the curve in a
point, $B$.

\subsection{Geometric formulation of a patch}
Consider \eqref{Eq:P_Patch} from the point of view of
projective geometry.
Given a patch
$\beta=\{\beta_\ba(x)\mid\ba\in\calA\}$ with $\calA\subset\R^d$, then 
 \begin{equation}\label{Eq:blending_param}
    x\ \longmapsto\ [\beta_\ba(x)\mid \ba\in\calA]
 \end{equation}
defines a map  $\beta\colon\Delta\to\R\P^\calA$ whose $\ba$th coordinate
is the blending function $\beta_\ba(x)$.

The map $\beta$ is unchanged if we
multiply all blending functions by the same positive function $\varphi(x)$
defined on $\Delta$.
Thus we
obtain the same map if we use instead the normalized blending functions.
The image $X_\beta$ of $\Delta$ under $\beta$ is a closed subset of
$\R\P^\calA$ (but not Zariski closed!), as $\Delta$ is compact and the map $\beta$ is
continuous.  
It is non-degenerate (does not lie in a hyperplane) exactly when the
patch is non-degenerate.

Suppose now that we have control points 
$b=\{\bb_\ba\in\R^\ell\mid \ba\in\calA\}$ and 
consider the composition of the map $\beta$ and the linear projection $\pi_b$
\[
  \Delta\ \xrightarrow{\ \beta\ }\ X_\beta\ \subset\ 
  \R\P^\calA\ \stackrel{\pi_b}{\longrat}\ \R\P^\ell\,.
\]
By our assumption on the positivity of the blending functions,
formula~\eqref{Eq:CoordsRd} implies that the image lies in the standard copy of $\R^\ell$, and 
there it is given by the formula~\eqref{Eq:P_Patch}.
Thus this image, $F(\Delta)$, is the image of $X_\beta$ 
under the projection $\pi_b$ given by the control points.

A reparametrization of the patch $\beta$ by a homeomorphism
$\psi\colon\Delta\to\Delta$ gives a different map
with the same image in $\R\P^\calA$
 \[
  \widetilde{\beta}\ \colon\ 
   \Delta\ \xrightarrow{\ \psi\ }\ 
   \Delta\ \xrightarrow{\ \beta\ }\  X_\beta\,.
 \]
Thus $X_\beta$ is an invariant of the patch modulo reparametrization.

If we take $\calA$ to be our set of control points where the element $\ba\in\calA$ is the
control point associated to the blending function $\beta_\ba$, then
the resulting linear projection $\pi_\calA$ is the 
\Blue{{\it tautological projection}}, 
 \begin{equation}\label{Eq:Taut_proj}
   \pi_\calA\ \colon\ \R\P^\calA\ \longrat\ \R\P^d\,.
 \end{equation}
Here, the coordinate points $(e_\ba\mid \ba\in\calA)$ are indexed by elements of $\calA$
and $\pi_\calA(e_\ba)=\ba\in\R\P^d$.
The tautological map $\tau$~\eqref{Eq:Tautological} is the composition
\[
  \Delta\ \xrightarrow{\ \beta\ }\ X_\beta\ \subset\ 
  \R\P^\calA\ \stackrel{\pi_\calA}{\longrat}\ \R\P^d\,.
\]
The patch is proper if this map is a homeomorphism onto its image
$\Delta\subset\R\P^d$, and the patch has linear precision if the
composition is the identity map, by Proposition~\ref{P:LP_criterion}.

With these definitions, we have the following identification of the reparametrization
having linear precision, which is a more precise version of Theorem~\ref{T:Lin_Prec}.

\begin{theorem}\label{Th:LP}
 Suppose that $X_{\beta}$ is a proper patch and $\pi_{\calA}\colon X_{\beta}\to\Delta$
 is the tautological projection restricted to $X_{\beta}$.
 Then the blending functions for $X_{\beta}$ that have linear precision are given by 
 the coordinates of the inverse of $\pi_{\calA}$. 
\end{theorem}

\begin{remark}\label{R:toric}
 This geometric perspective, where a patch is first an embedding
 of $\Delta$ into $\R\P^\calA$ followed by a linear projection, may have been introduced into
 geometric modeling in~\cite{dR91}.
 It is fundamental for Krasauskas's toric 
 B\'ezier patches~\cite{KR02}, and was reworked in the
 tutorial~\cite{So03}.
 There, an analog of Theorem~\ref{T:Lin_Prec} was formulated.
 We take this opportunity to correct an error in notation. 
 The tautological projection is an algebraic
 version of the moment map of the toric variety underlying 
 a toric patch, and not the actual moment map from symplectic geometry.\hfill\QED
\end{remark}

\subsection{Rational varieties}

We study the algebraic relaxation of our previous notions, replacing the real numbers by
the complex numbers so that we may use notions from algebraic geometry.
A rational patch has a parame\-trization by rational functions.
Clearing their denominators gives a new collection 
$\{\beta_\ba\mid\ba\in\calA\}$ of blending functions which are
polynomials. 
These polynomials define complex-valued functions on $\C^d$, and so the 
blending functions give a map
\[
  \beta_\C\ \colon\ \C^d\setminus B\ \longrightarrow\ \C\P^\calA\,,
\]
which is defined on the complement of the set 
$B$ where all the blending functions vanish.
This is called the \Blue{{\it base locus}} of the map $\beta_\C$.
This map $\beta_\C$ extends the map $\beta\colon\Delta\to X_\beta$.
We write $\beta$ for $\beta_\C$, $\P^\calA$ for $\C\P^\calA$, and in general
use the same notation for maps defined on complex algebraic varieties as for 
their restrictions to subsets of their real points.

As with linear projections, we write $\beta\colon\C^d\rat\P^\calA$ 
to indicate that $\beta$ is only defined on $\C^d\setminus B$.
Such a map between algebraic varieties that is given by polynomials and
defined on the complement
of an algebraic subset is called a \Blue{{\it rational map}}.

Let $Y_\beta$ be the Zariski closure of the image of $\beta$.
One important consequence of the patch being rational is that 
$X_\beta$ is a full-dimensional subset of the real points of the algebraic variety 
$Y_\beta$.
In particular, $X_\beta$ is defined locally in $\R\P^\calA$ by the
vanishing of some polynomials, and this remains true for any image
$F(\Delta)$ of a map $F$~\eqref{Eq:P_Patch} given by our original blending
functions and any choice of control points.
This is important in modeling, for these implicit equations are used
to compute intersections of patches.

Suppose that $X_\beta$ has rational linear precision and that $\beta=\{\beta_\ba\mid\ba\in\calA\}$
are rational blending functions with linear precision.
By Proposition~\ref{P:LP_criterion}, the composition
\[
  \Delta\ \xrightarrow{\ \beta\ }\ 
  X_\beta \ \stackrel{\pi_\calA}{\longrat}\ \Delta
\]
is the identity map.
Since $\Delta$ has full dimension in $\R^d$, it is Zariski dense in $\C^d$, and so 
the composition of rational functions
 \begin{equation}\label{Eq:comp}
  \C^d\ \stackrel{\beta}{\longrat}\ Y_\beta\ 
  \stackrel{\pi_\calA}{\longrat}\ \C^d
 \end{equation}
is also the identity map (where it is defined).

We introduce some terminology to describe this situation.
A rational map $\varphi\colon Y\rat Z$ between complex algebraic
varieties is a \Blue{{\it birational isomorphism}} if there is another map 
$\psi\colon Z\rat Y$ such that $\varphi\circ\psi$ and $\psi\circ\varphi$
are the identity maps where they are defined.  
In particular, $Y$ has a Zariski open subset $U$ and $Z$ a Zariski open
subset $V$ such that $\varphi|_U$ is one-to-one between $U$ and
$V$. 
The projection $\pi\colon C\rat L$ of Figure~\ref{F:cubic_Projection}
is not a birational isomorphism.
Indeed, if $x\in C$ is a point in $\pi^{-1}(p)$, then
$\overline{p,E}=\overline{x,E}$ and this plane meets 
the complement $C\setminus B$ of the base locus in $x$ and one other point.

Thus if the composition~\eqref{Eq:comp} is the identity map, then the map 
$\pi_\calA$ is a birational isomorphism from $Y_\beta\rat\C\P^d$.
We deduce the algebraic-geometric version of Proposition~\ref{P:LP_criterion}.

\begin{theorem}\label{Th:LP=birational}
 A patch $\beta=\{\beta_\ba\mid \ba\in\calA\}$ has rational linear precision if
 and only if the complexified tautological projection
\[
   \pi_\calA\ \colon\  Y_\beta\ \ \longrat\ \ \C\P^d
\]
 is a birational isomorphism.
\end{theorem}

Using well-known properties of birational projections, we use this 
characterization of linear precision to deduce 
necessary conditions for a patch to have linear precision.
A \Blue{{\it rational variety}} of dimension $d$ is one that is birational to
$\C^d$.

\begin{theorem}\label{Th:Necessary}
 If a patch $\beta=\{\beta_\ba\mid \ba\in\calA\}$ has rational linear precision, 
 then
 \begin{enumerate}
  \item $Y_\beta$ is a rational variety, 
  \item almost all codimension $d$ planes $L$ containing the center $E_\calA$ of the tautological 
        projection meet
        $Y_\beta$ in at most one point outside of $E_\calA$.
 \end{enumerate}
\end{theorem}

By condition (2), $Y_\beta$ has an exceptionally singular position with respect
to $E_\calA$. 
Typically, $Y_\beta$ does not meet a given codimension $d{+}1$ plane and its intersection
with a given codimension $d$ plane consists of $\deg(Y_\beta)$ points, counted with
multiplicity.
By Condition (2), not only does $Y_\beta$ meet $E_\calA$, but if $L$ has codimension $d$
and contains $E_\calA$, then most of the $\deg(Y_\beta)$ points in $Y_\beta\cap L$ 
lie in $E_\calA$. 
We will see this in Examples~\ref{Ex:WachH} and~\ref{Ex:Pentagon}.
\smallskip

\noindent{\it Proof.} 
 Statement (1) is immediate from Theorem~\ref{Th:LP=birational}.
 
 Since the birational isomorphism 
$\pi_\calA\colon Y_\beta\,\;\rat\;\,\C\P^d$  is the restriction of a linear
projection with center $E_\calA$, 
then for almost all points $y\in Y_\beta\setminus E_\calA$, we have 
\[
   \bigl(\overline{y,\,E_\calA}\bigr) \cap Y_\beta\ =\ 
    \{y\}\cup \bigl( E_\calA\cap Y_\beta\bigr)\,.
\]
 Since all codimension $d$ planes $L$ that contain $E_\calA$ have the form
 $(\overline{y,\,E_\calA})$ for some $y\not\in E_\calA$ and almost all meet $Y_\beta$,
 statement (2) follows.  
\hfill\QED

Consider the geometric situation of Figure~\ref{F:cubic_Projection}.
Suppose now that $E$ is the line tangent to $C$ at the point $B$.
Then $E\cap C$ is the point $B$, but counted with multiplicity 2.
A linear projection with center $E$ restricts to a birational isomorphism of $C$
with $\P^1$.

\subsection{Unirational varieties}

The necessary condition of Theorem~\ref{Th:Necessary}(1)
gives an important but subtle geometric restriction on patches that have
rational linear precision.

The map $\beta$ provides a parametrization of an open subset of $Y_\beta$
by an open subset of $\C^d$.
Such parametrized algebraic varieties are called 
\Blue{{\it unirational}}.
Unirational curves and surfaces are also rational, but these two notions
differ for varieties of dimension three and higher.
Clemens and Griffiths~\cite{CG72} showed that a smooth 
hypersurface of degree 3 in $\P^4$ is not rational (these were classically
known to be unirational, e.g. by Max Noether).

Thus parametric patches of dimension at least 3 will in general be unirational and not
rational. 
This does not occur for B\'ezier simploids or toric
patches.

\begin{example}\label{Ex:WachH}
 Consider again the Wachspress coordinates~\eqref{Eq:WHex}
 of Example~\ref{Ex:Wach}, which have linear precision.
 Let $Y_\beta\subset\P^5$ be the image of the blending functions~\eqref{Eq:WHex}.
 We study the base locus of the tautological projection $\pi_\calA$ of this patch.
 We first find equations for $Y_\beta$ as an algebraic subvariety of $\P^5$.

 If we divide
      $\beta_{\kzwei{-1}{-1}}\beta_{\kzwei{0}{1}}
       +2\beta_{\kzwei{0}{1}}\beta_{\kzwei{\smi 0}{-1}}
       +\beta_{\kzwei{1}{1}}\beta_{\kzwei{\smi 0}{-1}}$ 
 by the product of the linear forms defining the hexagon
 we get
\[
   (1-x)(1+y-x) + 2(1+x)(1-x) + (1+x)(1+x-y)
  \ =\ 4-2xy\,,
\]
 which is symmetric in $x$ and $y$,
 and so 
 \begin{equation}\label{Eq:beta_relation}
   \beta_{\kzwei{-1}{-1}}\beta_{\kzwei{0}{1}}
       +2\beta_{\kzwei{0}{1}}\beta_{\kzwei{\smi 0}{-1}}
       +\beta_{\kzwei{1}{1}}\beta_{\kzwei{\smi 0}{-1}}
   \ =\ 
   \beta_{\kzwei{-1}{-1}}\beta_{\kzwei{1}{0}}
        +2\beta_{\kzwei{1}{0}}\beta_{\kzwei{-1}{\smi 0}}
       +\beta_{\kzwei{1}{1}}\beta_{\kzwei{-1}{\smi 0}}\,.
 \end{equation}
 If we let $(y_\ba\mid\ba\in\calA)$ be natural coordinates for $\P^\calA$, 
then~\eqref{Eq:beta_relation} gives the quadratic polynomial 
which vanishes on $Y_\beta$,
\[
   y_{\kzwei{-1}{-1}}y_{\kzwei{0}{1}}
       +2y_{\kzwei{0}{1}}y_{\kzwei{\smi 0}{-1}}
       +y_{\kzwei{1}{1}}y_{\kzwei{\smi 0}{-1}}
   \ -\ 
   \bigl(y_{\kzwei{-1}{-1}}y_{\kzwei{1}{0}}
        +2y_{\kzwei{1}{0}}y_{\kzwei{-1}{\smi 0}}
       +y_{\kzwei{1}{1}}y_{\kzwei{-1}{\smi 0}}\bigr)\,.
\]
 Cyclically permuting the vertices of the hexagon gives two other
 quadratics that vanish on  $Y_\beta$, but these three sum to 0.
 There is another quadratic polynomial vanishing on $Y_\beta$,
\[
   y_{\kzwei{-1}{-1}}y_{\kzwei{0}{1}} + 
   y_{\kzwei{0}{1}}y_{\kzwei{1}{0}}+ 
   y_{\kzwei{1}{0}}y_{\kzwei{-1}{-1}}\ -\ 
  \bigl( y_{\kzwei{1}{1}}y_{\kzwei{\smi 0}{-1}} + 
         y_{\kzwei{\smi 0}{-1}}y_{\kzwei{-1}{\smi 0}}+ 
         y_{\kzwei{-1}{\smi 0}}y_{\kzwei{1}{1}}\bigr)\ .
\]

 There is an additional cubic relation among the blending functions, 
 which gives a cubic polynomial that vanishes on $Y_\beta$
 \begin{equation}\label{Eq:cubic}
   y_{\kzwei{-1}{-1}}y_{\kzwei{0}{1}}y_{\kzwei{1}{0}}\ -\ 
   y_{\kzwei{1}{1}}y_{\kzwei{\smi 0}{-1}}y_{\kzwei{-1}{\smi 0}}\,.
 \end{equation}
 These relations, three independent quadratics and one cubic, define $Y_\beta$ as an
 algebraic subset of $\P^\calA=\P^5$.
 The quadratic equations define $Y_\beta$, together with the
 2-dimensional linear space cut out by the three linear equations.
 \begin{equation}\label{Eq:lin_forms}
 \begin{array}{rcl}
  0&=&
   y_{\kzwei{-1}{-1}} +  y_{\kzwei{0}{1}} + y_{\kzwei{1}{0}} + 
   y_{\kzwei{1}{1}} + y_{\kzwei{\smi 0}{-1}} +y_{\kzwei{-1}{\smi 0}}\\
  0&=&\rule{0pt}{13pt}
   y_{\kzwei{1}{1}} -y_{\kzwei{-1}{-1}}+y_{\kzwei{1}{0}}-y_{\kzwei{-1}{\smi 0}}\\
  0&=&\rule{0pt}{13pt}
  y_{\kzwei{1}{1}} -y_{\kzwei{-1}{-1}}+y_{\kzwei{0}{1}}-y_{\kzwei{\smi 0}{-1}}
 \end{array}
\end{equation}

 Thus $Y_\beta$ has degree $7=2^3-1$.
 A general codimension $2$ plane $L$ containing the center $E_\calA$ of the tautological projection
 $\pi_\calA$ 
 will meet $Y_\beta$ in at most 1 point outside of $E_\calA$, by Theorem~\ref{Th:Necessary}(2).
 Since $L\cap Y_\beta$ has degree at least 7, 
the other 6 points must lie in the base locus 
$E_\calA\cap Y_\beta$, which is in fact 
a reducible cubic plane curve.

 The three linear forms~\eqref{Eq:lin_forms} also define the tautological 
 projection, $\pi_\calA$, so this linear subspace is the center of projection
 $E_\calA$, and the base locus $B$ is defined in $E_\calA$ by the cubic~\eqref{Eq:cubic}
 defining $Y_\beta$.
 Let us parametrize the center $E_\calA$ as follows.
 For $[a,b,c]\in\P^2$, set
\[
  \begin{array}{rclcrcl}
    y_{\kzwei{-1}{-1}}   &:=& a+b  && y_{\kzwei{1}{1}}      &:=& -a+b\\
    y_{\kzwei{0}{1}}     &:=& a-c  && y_{\kzwei{\smi 0}{-1}}&:=& -a-c\\
    y_{\kzwei{1}{0}}     &:=& a-b+c&& y_{\kzwei{-1}{\smi 0}}&:=& -a-b+c
  \end{array}
\]
 Then we have 
\[
   y_{\kzwei{-1}{-1}}y_{\kzwei{0}{1}}y_{\kzwei{1}{0}}\ -\ 
   y_{\kzwei{1}{1}}y_{\kzwei{\smi 0}{-1}}y_{\kzwei{-1}{\smi 0}}\ =\ 
    2a(a^2-b^2+bc-c^2)\, ,
\]
 the product of a linear and a quadratic form.
 \hfill\QED
\end{example}

%
%
%
\section{Toric Patches}

Toric patches, which were introduced by Krasauskas~\cite{KR02}, are a
class of patches based on certain 
special algebraic varieties, called toric varieties.
For a basic reference on toric varieties as parametrized varieties, see the
book~\cite{GBCP}, particularly Chapters 4 and 13, and the articles~\cite{Cox03,So03},
which are introductions aimed at people in geometric modeling.
We first define toric patches, give some examples, reinterpret
Theorem~\ref{Th:LP=birational} for toric patches, and then state the main open problem
about linear precision for toric patches. 

\subsection{Toric Patches}
We regard elements $\ba=(a_1,\dotsc,a_d)\in\Z^d$ as exponent vectors for 
monomials in the variables $x=(x_1,\dotsc,x_d)$
\[
   x^\ba\ :=\ x_1^{a_1} x_2^{a_2} \dotsb x_d^{a_d}\,.
\]
Let $\calA\subset\Z^d$ be a finite set of integer vectors whose convex hull $\Delta$ is a 
polytope of dimension $d$ in $\R^d$.
This implies that $\calA$ affinely spans $\R^d$.
Suppose that $w=\{w_\ba\in\R_>\mid \ba\in\calA\}$ is a set of positive 
weights indexed by $\calA$.
These data $(\calA,w)$ define a map  
$\varphi_{\calA,w}\colon(\C^\times)^d\to\P^\calA$,
 \begin{equation}\label{Eq:monom_parametrize}
    \varphi_{\calA,w}\ \colon\ x \ \longmapsto
  \ [\,w_\ba x^\ba\mid \ba\in\calA\,]\,.
 \end{equation} 
Let $Y_{\calA,w}$ be the closure of the image of
$(\C^\times)^d$ under the map $\varphi_{\calA,w}$.
When the weights $w_\ba$ are equal, this is the toric variety
parametrized by the monomials in $\calA$.
In general, it is the translate of that toric variety by the element 
$w$ of the positive real torus $\R_>^\calA$, which acts on $\P^\calA$ by
scaling the coordinates.

Let \Blue{$X_{\calA,w}}\subset\R\P^\calA$ be the closure $\varphi_{\calA,w}(\R_>^d)$.
This is the non-negative part of the translated toric variety 
$Y_{\calA,w}$.
That is, if $\Blue{\R\P^\calA_\geq}\subset \P^\calA$ is the set of points having
non-negative (real) homogeneous coordinates, then 
$X_{\calA,w}=Y_{\calA,w}\cap\R\P^\calA_\geq$.
It is homeomorphic to the convex hull $\Delta$ of the 
vectors in $\calA$, as a manifold with corners~\cite[\S~4]{Fu93}.

The set $\calA$ is \Blue{{\it primitive}} if it affinely spans $\Z^d$.
We may assume that $\calA$ is primitive, for there is always some primitive set
$\calA'\subset\Z^d$ with $Y_{\calA,w}=Y_{\calA',w}$.
First note that translating every vector in $\calA$ by a fixed
vector $\ba'$ multiplies each coordinate of the map $\varphi_{\calA,w}$ by the
monomial $x^{\ba'}$, which does not change $\varphi_{\calA,w}$ as a  map to projective space.
Thus we may assume that $\calA$ contains the zero vector.
Since $\calA$ affinely spans $\R^d$, it generates a subgroup $\Z\calA$ of $\Z^d$ which is
isomorphic to $\Z^d$.
Let $\bb_1,\dotsc,\bb_d\in\Z\calA$ be the basis corresponding to the basis of $\Z^d$
under an isomorphism $\Z\calA\xrightarrow{\,\sim\,}\Z^d$.
If we let $\calA'\subset\Z^d$ be the image of $\calA$ under this isomorphism 
and $\calB:=\{\bb_1,\dotsc,\bb_d\}$, 
then the map $\varphi_{\calA,w}$ factors
\[
   (\C^\times)^d\ \xrightarrow{\ \chi_\calB\ }
   (\C^\times)^d\ \xrightarrow{\ \varphi_{\calA',w}\ }\P^\calA\,,
\]
where $\chi_\calB$ is the surjective map that sends $x$ to 
$(x^{\bb_1},\dotsc,x^{\bb_d})$.
Note that $\varphi_{\calA', w}$ is injective but that $\chi_{\calB}$ has 
fibers of size $|\Z^{d}/\Z\calA|$.
Thus $\varphi_{\calA,w}$ is not injective unless $\calA$ is primitive.

\begin{definition}\label{DE:TP}
 A \Blue{{\it toric patch}} of \Blue{{\it shape $(\calA,w)$}} is any patch 
 $\{\beta_\ba\mid\ba\in\calA\}$
 such that the image, $X_\beta$, of 
 $\Delta$ under the map~\eqref{Eq:blending_param} given by the blending
 functions is equal to the non-negative part $X_{\calA,w}$ of the
 translated toric variety $Y_{\calA,w}$.\hfill\QED
\end{definition}

\begin{remark}\label{Rem:weights}
  Our geometric perspective that the fundamental object is the image $X_\beta$ requires us
  to absorb the weights $w$ into our definition and distinguishes patches with different
  choices of weights.\hfill\QED
\end{remark}

\begin{example}\label{Ex:TBC}
 Let $\calA:=\{0,1,\dotsc,n\}\subset\Z$ and let $w\in\R_>^{n+1}$ be any
 set of weights.
 Then the map $\varphi_{\calA,w}$ is
\[
  x\ \longmapsto\ [w_0,\, w_1x,\, w_2x^2,\,\dotsc,\, w_nx^n]\,.
\]
 The map $x\mapsto \frac{x}{n-x}$ sends $\Delta=[0,n)$ to $[0,\infty)$.
 Precomposing $\varphi_{\calA,w}$ with this map gives a
 toric patch of shape $(\calA,w)$.
\newcommand{\xnx}{\bigl({\textstyle\frac{x}{n-x}}\bigr)}
 \begin{eqnarray*}
  [0,n]\ \ni\ x&\mapsto&\Bigl[w_0,\,w_1\xnx,\,w_2\xnx^2,\,\dotsc,\,w_n\xnx^n\Bigr]\\
   &=&[w_0(n-x)^n,\, w_1x(n-x)^{n-1},\,w_2x^2(n-x)^{n-2},\,
   \dotsc,w_nx^n]\,.
 \end{eqnarray*}
Observe that if we choose weights $w_i=\binom{n}{i}$, substitute 
$ny$ for $x$, and remove the common factors of $n^n$, then we obtain the
blending functions for the B\'ezier curve of Example~\ref{Ex:BC}.
Note that replacing $x$ by $ny$ also replaces
$i\in\calA$ by $\frac{i}{n}$. 
\hfill\QED
\end{example}

\begin{example}\label{Ex:TBP}
 Krasauskas~\cite{KR02} generalized the classical B\'ezier
 parametrization to any polytope with integer vertices.
 A polytope $\Delta$ is defined by its facet inequalities
\[
  \Delta\ =\ \{x\in\R^d\mid h_i(x)\geq 0\,, i=1,\dotsc,N\}\,.
\]
 Here, $\Delta$ has $N$ facets and for each $i=1,\dotsc,N$, 
 $h_i(x)=\langle \bv_i,x\rangle + c_i$ is the linear function defining
 the $i$th facet, where $\bv_i\in\Z^d$ is the (inward oriented) primitive vector normal to 
 the facet and $c_i\in\Z$. 
 (Compare these to the functions in Example~\ref{Ex:Wach}.)

 Let $\calA\subset\Delta\cap\Z^d$ be any subset of the integer points of
 $\Delta$ that includes its vertices.
 Let $w=\{w_\ba\mid \ba\in\calA\}\subset \R_>$ be a collection of
 positive weights.
 For every $\ba\in\calA$, Krasauskas defined \Blue{{\it toric B\'ezier functions}}
 \begin{equation}\label{Eq:toric-Bezier}
   \beta_\ba(x)\ :=\
    w_\ba h_1(x)^{h_1(\ba)}h_2(x)^{h_2(\ba)}\dotsb h_N(x)^{h_N(\ba)}\,.
 \end{equation}
 For $x\in\Delta$, these parametrize a patch $X_\beta$.

 Krasauskas defined a toric patch to be any reparametrization of such a
 patch.
 Our definition of a toric patch (Definition~\ref{DE:TP}) agrees with Krasauskas's,
 when $\calA=\Delta\cap\Z^d$ consists of all the integer points in a polytope:
 Observe that the map $\beta$ in~\eqref{Eq:toric-Bezier} 
 is the composition of a map $H:\Delta\to\R^{N}$ given by
 $x\mapsto(h_{1}(x),\dots,h_{N}(x))$  
 with a rational map $\varphi: \C^{N} \rat\ \R\P^{\calA}$ given by 
 \[ 
    \varphi\ \colon\ (u_{1},\dots,u_{N})\ \longmapsto\ 
      [w_{\ba}u_{1}^{h_{1}(\ba)}\dotsm u_{N}^{h_{N}(\ba)} \mid \ba \in \calA]\,. 
 \] 
 This map $\varphi$ factors through the map $\varphi_{\calA,w}$.
 Indeed, define a map
 \begin{eqnarray*}
  f_\Delta\ \colon\ (\C^\times)^N&\longrightarrow& (\C^\times)^n\\
   u=(u_1,\dotsc,u_N)&\longmapsto&t=(t_1,\dotsc,t_n)\quad\mbox{where}\quad
    t_j:=\prod_{i=1}^N u_i^{\langle \bv_i,\be_j\rangle}
  \end{eqnarray*}
 Then
\[
   \prod_{i=1}^N u_i^{h_i(\ba)}\ =\ \prod_{i=1}^N u_i^{c_i}\cdot
     \prod_{i=1}^N u_i^{\langle \bv_i,\ba\rangle}
      \ =\ u^c\cdot t^{\ba}\,.
\]
 And so $\varphi(u)=u^c\cdot\varphi_{\calA,w}(f_\Delta(u))$.
 This shows that $X_\beta\subset X_{\calA,w}$, and Krasauskas~\cite[\S 4]{Kr06} shows 
 that this inclusion is an equality.
\hfill\QED
\end{example}

\begin{example}\label{Ex:CBP}
 A \Blue{{\it standard $d$-simplex of degree $n$}} is the simplex in 
$\R^d$,
\[
   n\Delta_d\ :=\ 
   \{x\in \R^d\mid x_1+\dotsb+x_d\leq n\mbox{\ and\ }x_i\geq 0\}\,.
\]
 For a point $\ba=(a_1,\dotsc,a_d)$ in $n\Delta_d\cap\Z^d$ set $|\ba|:=a_1+\dotsb+a_d$
 and let $w_\ba:=\binom{n}{\ba}=\frac{n!}{a_1!\dotsb a_d! (n-|\ba|)!}$ be the multinomial  
 coefficient.
 This gives a system of weights for $\calA= n\Delta\cap\Z^d$.
 Then Krasauskas's toric B\'ezier
 patch for $n\Delta_d$ has blending functions
\[
   \beta_\ba\ :=\ {\textstyle \binom{n}{\ba}}
   (n-{\textstyle \sum_i x_i})^{n-\sum_i a_i}\,
    \prod_{i=1}^d x_i^{a_i}\,.
\]
 If we substitute $ny_i$ for $x_i$ and then remove
 the common factors of $n^n$, we recover the classical Bernstein
 polynomials~\cite[\S 4.1]{Fa97}. \hfill\QED
\end{example}

\begin{example}
 The \Blue{{\it  B\'ezier simploids}}~\cite{dRGHM} are toric patches based on 
 products of B\'ezier simplices. 
 Suppose that $d_1,\dotsc,d_m$ and $n_1,\dotsc,n_m$ are positive integers.
 Set $d:=d_1+\dotsb+d_m$.
 Write $\R^d$ as a direct sum $\R^{d_1}\oplus\dotsb\oplus\R^{d_m}$ and suppose that for each
 $i=1,\dotsc,m$, the scaled simplex $n_i\Delta_{d_i}$ lies in the $i$th summand $\R^{d_i}$.
 Consider the Minkowski sum
\[
  \Delta\ :=\ n_1\Delta_{d_1} + n_2\Delta_{d_2} + \dotsb + n_m \Delta_{d_m}\ \subset\ \R^d\,,
\]
 which is just the product of the simplices $n_i\Delta_{d_i}$.

 If $\bx_i:=(x_{i,1},\dotsc,x_{i,d_i})$ are coordinates on $\R^{d_i}$, then 
\[
   \Delta\ =\ \{ \bx\ =\ (\bx_1,\dotsc,\bx_m) \mid 
     x_{i,j}\geq 0\ \quad\mbox{and}\quad\ 
     x_{i,1}+\dotsb+x_{i,d_i}\leq n_i\quad i=1,\dotsc,m\}\,.
\]
 Let $\calA:=\Delta\cap\Z^d$ be the set of integer points  of $\Delta$.
 Given an integer point $\ba=(\ba_1,\dotsc,\ba_m)\in\calA$,
 set $w_{\ba}:=\prod_{i=1}^m \binom{n_i}{\ba_i}$, 
 the product of multinomial coefficients.
 This gives a system of weights, and let $Y_{\calA,w}$ be the resulting toric patch.
 Krasauskas's toric patch of shape $\Delta$ has blending
 functions  
\[
   \beta_{\ba}(\bx)\ :=\ \prod_{i=1}^m \beta_{\ba_i}(\bx)\,,
\]
 where $\beta_{\ba_i}(\bx_i)$ for $\bx_i\in\R^{d_i}$ are the Bernstein polynomials.
 The resulting patch is called a \Blue{{\it B\'ezier simploid}} in~\cite{dRGHM}.
\hfill\QED
\end{example}

We state our main result about linear precision for toric patches,
which is a useful reformulation of linear precision for toric patches.

\begin{theorem}\label{Th:toric_diff}
  Let $\calA\subset\Z^d$ be a primitive collection of exponent vectors, 
  $w\in\R_>^{\calA}$ be a system of weights, and define the Laurent polynomial
\[
   f\ =\ \Blue{f_{\calA,w}}\ :=\ \sum_{a\in\calA} w_{a} x^a\,.
\]
 A toric patch of shape $(\calA,w)$ has rational linear precision if and only if the
 rational function  $\psi_{\calA,w}\colon\C^d\rat\C^d$
 defined by
 \begin{equation}\label{Eq:toric_differential}
  \frac{1}{f}\Bigl( x_1 \frac{\partial}{\partial x_1}f,\ 
                    x_2 \frac{\partial}{\partial x_2}f,\ \dotsc,\ 
                    x_d \frac{\partial}{\partial x_d}f \Bigr)
 \end{equation} 
 is a birational isomorphism.
\end{theorem}

\begin{remark}
The map~\eqref{Eq:toric_differential} has an interesting reformulation in terms of toric
derivatives. 
The $i$th \Blue{{\it toric derivative}} of a Laurent polynomial $f$ is
$x_i\frac{\partial}{\partial x_i}f$. 
The toric differential $D_{\mbox{\scriptsize\rm torus}}f$ is the vector whose components are the 
toric derivatives of $f$.
Thus the map $\psi_{\calA,w}$ is the logarithmic toric differential
$D_{\mbox{\scriptsize\rm torus}} \log f$.

In the proof below, we show that the map $\psi_{\calA,w}$ comes from a map $\C^d\;\rat\;\C\P^d$,
 \begin{equation}\label{Eq:MAP}
   x\ \longmapsto\ 
   \Bigr[f\;:\;  x_1 \frac{\partial}{\partial x_1}f\;:\;
                    x_2 \frac{\partial}{\partial x_2}f\;:\ \dotsb\ :\;
                    x_d \frac{\partial}{\partial x_d}f \bigr]
 \end{equation}
Let $F_{\calA,w}$ be the homogenization of $f$ with respect to a new variable $x_0$. 
Then 
\[
   \deg(F_{\calA,w})F_{\calA,w}\ =\ \sum_{i=0}^d x_i\frac{\partial}{\partial x_i}F_{\calA,w}\,.
\]
Thus, after homogenizing and a linear change of coordinates, the map~\eqref{Eq:MAP} 
is defined by the formula
\[
  \Blue{D_{\mbox{\scriptsize\rm torus}}F_{\calA,w}}\ =\ 
   \Bigr[x_0 \frac{\partial}{\partial x_0}F_{\calA,w}\;:\;
         x_1 \frac{\partial}{\partial x_1}F_{\calA,w}\;:\ \dotsb\ :\
         x_d \frac{\partial}{\partial x_d}F_{\calA,w} \bigr]
  \makebox[.1cm][l]{\hspace{2cm}\QED}
\]
\end{remark}

This last form of the map~\eqref{Eq:MAP} gives an appealing reformulation
of Theorem~\ref{Th:toric_diff} which we feel is the most useful for  
further analysis of linear precision for toric patches.

\begin{corollary}\label{Cor:Reformulation}
 The toric patch of shape $(\calA,w)$ has linear precision if and only if the map
 $D_{\mbox{\scriptsize\rm torus}}F_{\calA,w}\colon\C\P^{d}\rat\C\P^{d}$
 is a birational isomorphism.
\end{corollary}

\noindent{\it Proof of Theorem~$\ref{Th:toric_diff}$.}
 By Theorem~\ref{Th:LP=birational}, a toric patch of shape $(\calA,w)$
 has linear precision if and only if the complexified tautological projection
\[
   \pi_\calA\ \colon\ Y_{\calA,w}\ \ \longrat\ \ \C\P^d
\]
 is a birational isomorphism.
 Precomposing this with the defining parametrization 
 of $Y_{\calA,w}$,  $\varphi_{\calA_w}\colon(\C^\times)^d\hookrightarrow Y_{\calA,w}$
 (which is a rational map $\C^d\rat Y_{\calA,w}$), gives a map
\begin{equation}\label{Eq:composite}
  \begin{array}{rcccl}
    \C^d\ &\ \stackrel{\varphi_{\calA,w}}{\longrat}\ &Y_{\calA,w}&
     \ \stackrel{\pi_{\calA}}{\longrat}\ &\C\P^d \\
     x\ &\ \longmapsto & [w_\ba x^\ba\mid \ba\in\calA]& \longmapsto & \rule{0pt}{15pt}
   {\displaystyle \sum_{\ba\in\calA} w_\ba x^\ba (1,\ba)}\,.
  \end{array}
\end{equation}

 The initial (0th) coordinate of this map is
\[
   \sum_{\ba\in\calA} w_\ba x^\ba
\]
 which is the polynomial $f_{\calA,w}$.
 For $i>0$ the $i$th coordinate of this composition is
\[
  \sum_{\ba\in\calA} a_i\cdot w_\ba x^\ba\ =\  x_i\frac{\partial}{\partial x_i}f_{\calA,w}\,.
\]
 If we divide the composite map~\eqref{Eq:composite} by its 0th coordinate, $f_{\calA,w}$, 
 we obtain the map $\psi_{\calA,w}$.
 This map is birational if and only if $\pi_\calA$ is birational, as $\varphi_{\calA,w}$ is
 birational. 
\hfill\QED

\begin{example}
 Let $\calA:=\{0,1,\dotsc,n\}\subset\Z$ as in Example \ref{Ex:TBC}.
 Let $w\in\R_>^{n+1}$ be any set of weights.
 Then the polynomial $f = f_{\calA,w}$ is
\[
  f(x) \ = \ w_0 + w_1x + \dotsm + w_nx^n\,.
\]

Suppose that this toric patch has rational linear precision.
By Theorem~\ref{Th:toric_diff}, the logarithmic toric differential 
$xf'(x)/f(x)$  is a birational isomorphism. 
Thus the fraction $xf'(x)/f(x)$ reduces to a quotient of linear
polynomials and so $f(x)$ is necessarily a pure power, say $a(x+\alpha)^n$ 
with $\alpha\neq 0$.

We conclude that the weights are $w_i=\alpha^{n-i}\binom{n}{i}$, which are not quite the
weights in the B\'ezier curve~\ref{Ex:BC}.
If we rescale $x$ (which does not change $Y_{\calA,w}$), setting $y=x/\alpha$, then the
weights become $|\alpha|^n\binom{n}{i}$. 
Removing the common factor of $|\alpha|^n$, reveals that this 
toric patch is the classical B\'ezier curve of Example~\ref{Ex:BC}. 
\hfill\QED
\end{example}

We close with two problems.\medskip

\noindent{\bf Problem 1.}
Classify the toric patches of dimension $d$ that have linear precision.\medskip

The analysis of this section was recently used to classify toric surface patches ($d=2$) 
having linear precision~\cite{BRS}.
It remains an open problem to understand how to tune a toric patch (moving the points
$\calA$) to achieve linear precision.\medskip 

\noindent{\bf Problem 2.}
Classify the toric patches of dimension $d$ that may be tuned to have linear
precision.\medskip

This is a strictly larger class of patches.
For example, the classification of~\cite{BRS} shows if a toric surface patch has linear
precision, then it is triangular or quadrilateral,
but it is possible to tune a pentagonal patch to achieve linear precision.

\begin{example}\label{Ex:Pentagon}
 We discuss Krasauskas's version~\cite{Kr06} of  Kar{\v{c}}iauskas's pentagonal
 $M$-patch~\cite{Ka03}, which has linear precision.
 Let $\calA$ be the integer points in the 
 pentagon~\includegraphics{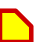} shown below, 
 and let $\{w_\ba\mid \ba\in\calA\}$ be the indicated weights.
\[
  \begin{picture}(125,120)(-10,-18)
   \put(17,11){\includegraphics{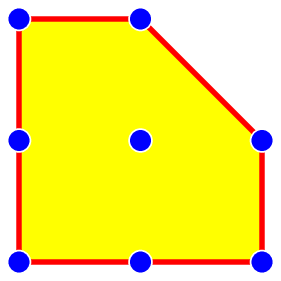}}
   \put(-10,89){$(0,2)$} \put(40,89){$(1,2)$}
   \put(-12,44){$(0,1)$} \put(40,55){$(1,1)$}\put(92,44){$(2,1)$}
   \put(-10,-1){$(0,0)$} \put(40,-1){$(1,0)$}\put(87, -1){$(2,0)$}
   \put(10,-18){Lattice points $\calA$}
  \end{picture}
  \qquad\qquad\quad
  \begin{picture}(110,120)(-10,-18)
   \put(4,11){\includegraphics{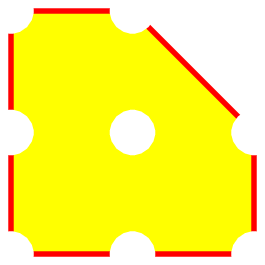}}
   \put(2,78){$2$}  \put(37,78){$2$} 
   \put(2,43){$5$}  \put(37,43){$7$} \put(72,43){$2$}
   \put(2, 8){$3$}  \put(37, 8){$5$} \put(72, 8){$2$}
   \put(-10,-18){weights $\{w_\ba\mid\ba\in\calA\}$}
  \end{picture}
\]
 These define a toric patch $X_{\calA,w}$ and translated toric variety $Y_{\calA,w}$.
 This subvariety of $\P^7$ has degree 7 and is defined by 14 quadratic equations, all of
 the form  
\[
   y_\ba y_\bb-y_\bc y_\bd\qquad
   \mbox{where}\quad \ba,\bb,\bc,\bd\in\calA
   \quad\mbox{with}\quad \ba+\bb=\bc+\bd\,.
\]

 This has a monomial parametrization 
\[
    \varphi_{\calA,w}\ \colon\ (s,t)\ \longmapsto\ 
     [s^{a_1}t^{a_2}\,\mid\, \ba=(a_1,a_2)\in\calA]\,.
\]
 Its composition with the tautological projection $\pi_\calA$ gives a 
 map $g\colon\P^2\rat\P^2$ of degree 4, with 
 3 base points at $(-1,-1)$, $(-1,-\frac{1}{2})$, and $(-\frac{1}{2},-1)$

 We tune this patch (or rather the projection $\pi_\calA$) by moving the
  non-extreme points of $\calA$ within the pentagon.
 Let $\calB$ consist of the vertices of the pentagon, together with the 
 three non-vertices
\[
  \begin{picture}(125,120)(-10,-20)
   \put(17,11){\includegraphics{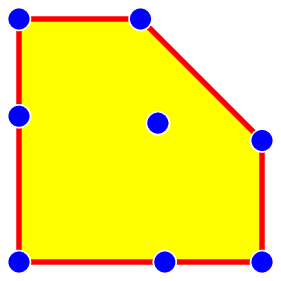}}
   \put(-10,89){$(0,2)$} \put(40,89){$(1,2)$}
   \put(-14,51){$(0,\frac{6}{5})$} 
     \put(30,37){$(\frac{8}{7},\frac{8}{7})$}\put(92,45){$(2,1)$}
   \put(-10,-2){$(0,0)$} \put(47,-2){$(\frac{6}{5},0)$}\put(85,-2){$(2,0)$}
   \put(10,-20){Lattice points $\calB$}
  \end{picture}
\]
 The composition of $\varphi_{\calA,w}$ with this new tautological projection $\pi_\calB$ 
 gives the rational map 
\[
   (s,t)\ \longmapsto\ 
    \Bigl( \frac{2s(2s+t+3)}{(s+1)(2s+2t+3)}\,,\ 
           \frac{2t(s+2t+3)}{(t+1)(2s+2t+3)} \Bigr)\,.
\]
 This is given by quadrics and has 3 base points $(-1,-1)$, $(0,-\frac{3}{2})$, 
 and $(-\frac{3}{2},0)$,
 and so it has degree $1=2\cdot 2-3$, and is therefore a birational isomorphism.
 By Theorem~\ref{Th:LP=birational}, the tuned toric patch 
 (with the lattice points $\calA$ replaced by $\calB$) has linear precision.

Thus there is a rational parametrization
$\beta\colon\includegraphics[height=10pt]{figures/pent.eps}\to X_{\calA,w}$ 
whose composition with the tautological projection $\pi_\calB$ is the identity map.
Indeed, Krasauskas~\cite{Kr06} gives the following modification of 
the toric B\'ezier functions and shows that the map 
$\pi_\calB\circ\beta$ is the identity.
  \begin{equation}\label{Eq:NTBF}
   \begin{array}{rcl}
   \beta_{\kzwei{0}{0}} &:=& 3(2-s)^2(3-s-t)^3(2-t)^2 \\\rule{0pt}{14pt}
   \beta_{\kzwei{1}{0}} &:=& 5(2-s)(3-s-t)^2(2-t)^2s(3-s-\tfrac{t}{2}) \\\rule{0pt}{14pt}
   \beta_{\kzwei{0}{1}} &:=& 5(3-\tfrac{s}{2}-t)t(2-s)^2(3-s-t)^2(2-t) \\\rule{0pt}{14pt}
   \beta_{\kzwei{1}{1}} &:=& 7(3-\tfrac{s}{2}-t)t(2-t)(3-s-t)(2-s)s(3-s-\tfrac{t}{2})
           \\\rule{0pt}{14pt}
   \beta_{\kzwei{2}{0}} &:=& 2(3-s-t)(2-t)^2s^2(3-s-\tfrac{t}{2})^2 \\\rule{0pt}{14pt}
   \beta_{\kzwei{0}{2}} &:=& 2(3-\tfrac{s}{2}-t)^2t^2(2-s)^2s(3-s-t) \\\rule{0pt}{14pt}
   \beta_{\kzwei{2}{1}} &:=& 2(3-\tfrac{s}{2}-t)t(2-t)s^2(3-s-\tfrac{t}{2})^2 \\\rule{0pt}{14pt}
   \beta_{\kzwei{1}{2}} &:=& 2(3-\tfrac{s}{2}-t)^2t^2(2-s)s(3-s-\tfrac{t}{2}) 
 \end{array}
\end{equation}
This is a modification of the toric B\'ezier functions for the pentagon.
If, in the definition~\eqref{Eq:toric-Bezier}, we  replace the form $s$ 
by $s(3-s-\tfrac{t}{2})$ and the form $t$ by
$(3-\tfrac{s}{2}-t)t$,  then we get these modified blending
functions~\eqref{Eq:NTBF}. \smallskip

 Let us consider the geometry of linear precision for this tuned patch in the spirit of 
 Section~\ref{S:Geometry}.
 The projection $\pi_\calB$ is defined by the three linear forms.
 \begin{eqnarray*}
   X&:=& \tfrac{6}{5}y_{\kzwei{1}{0}}\;+\;2y_{\kzwei{2}{0}}\;+\;
         \tfrac{8}{7}y_{\kzwei{1}{1}}\;+\;2y_{\kzwei{2}{1}}\;+\;y_{\kzwei{1}{2}}\\
   Y&:=& \tfrac{6}{5}y_{\kzwei{0}{1}}\;+\;2y_{\kzwei{0}{2}}\;+\;
         \tfrac{8}{7}y_{\kzwei{1}{1}}\;+\;2y_{\kzwei{1}{2}}\;+\;y_{\kzwei{2}{1}}\\
   Z&:=&y_{\kzwei{0}{0}}\;+\;y_{\kzwei{1}{0}}\;+\;y_{\kzwei{0}{1}}\;+\;y_{\kzwei{1}{1}}\;+\;
        y_{\kzwei{2}{0}}\;+\;y_{\kzwei{0}{2}}\;+\;y_{\kzwei{2}{1}}\;+\;y_{\kzwei{1}{2}}
 \end{eqnarray*}
Its center $E_\calB$ is defined by the vanishing of these three forms.
A general linear subspace $L$ of codimension 2 containing $E_\calB$ is defined by
 equations of the form 
\[
   X\ =\ x Z\qquad\mbox{and}\qquad Y\ =\ yZ\,.
\]
 This subspace meets $Y_{\calA,w}$ in four points.
 One point lies outside of $E_\calB$ while three points lie in the center $E_\calB$.
 We give the three points.
\[
  \begin{picture}(110,120)(0,-18)
   \put(4,11){\includegraphics{figures/pent_weights.eps}}
   \put(2,78){$3$}  \put(37,78){$0$} 
   \put(-7,43){$-5$}  \put(37,43){$0$} \put(72,43){$0$}
   \put(2, 8){$2$}  \put(37, 8){$0$} \put(72, 8){$0$}
   \put(4,-18){multiplicity 1}
  \end{picture}
  \qquad
  \begin{picture}(110,120)(0,-18)
   \put(4,11){\includegraphics{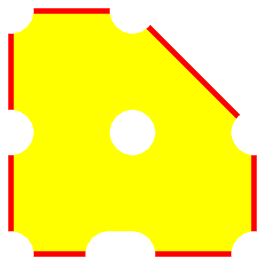}}
   \put(2,78){$0$}  \put(37,78){$0$} 
   \put(2,43){$0$}  \put(37,43){$0$} \put(72,43){$0$}
   \put(2, 8){$2$}  \put(28, 8){$-5$} \put(72, 8){$3$}
   \put(4,-18){multiplicity 1}
  \end{picture}
  \qquad
  \begin{picture}(110,120)(0,-18)
   \put(4,11){\includegraphics{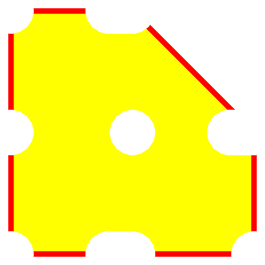}}
   \put(2,78){$2$}  \put(28,78){$-2$} 
   \put(-7,43){$-5$}  \put(37,43){$7$} \put(63,43){$-2$}
   \put(2, 8){$3$}  \put(28, 8){$-5$} \put(72, 8){$2$}
   \put(4,-18){multiplicity 4}
  \end{picture}
\]
 We write the coordinate $y_{\ba}$ of a point $y\in\R\P^{\calA}$ at the
 position of the point $\ba$.\hfill\QED
\end{example}

%
%
%
\section{Iterative proportional fitting for toric patches}

The toric patch $X_{\calA,w}$ from geometric modeling appears naturally in algebraic
statistics in the form of a toric model, which leads to a dictionary between 
the subjects.
We show that a toric patch has rational linear precision if and only if it has maximum
likelihood degree 1 as a statistical model. 
As a consequence, we present a new family
of toric patches with rational linear precision. 
Finally, we propose the iterative proportional fitting algorithm from 
statistics as a tool to compute the unique reparametrization of a toric patch having
linear precision.

In algebraic statistics, the image of $\R^{d}_{>}$ under the map
$\varphi_{\calA,w}$~\eqref{Eq:monom_parametrize}  
is known as a \Blue{{\it toric model}}~\cite[\S 1.2]{PS05}.
They are more commonly called \Blue{{\it log-linear models}}, as the logarithms of the 
coordinates of $\varphi_{\calA,w}$ are linear functions in the logarithms of the 
coordinates of $\R^{d}_{>}$, or \Blue{{\it discrete exponential families}}, as 
the coordinates of  $\varphi_{\calA,w}$ are exponentials of the coordinates of $\R^d$, which
are themselves logarithms of the coordinates of $\R^{d}_{>}$.

We identify the non-negative orthant $\R\P^\calA_{\geq}$  with the 
probability simplex
\[
   \Blue{\Delta^\calA}\ :=\ 
   \{y\in\R^\calA_{\geq0}\mid |y|:=\sum_{a\in\calA} y_a=1\}\,,
\]
and so we may also regard  $X_{\calA,w}$ as a subvariety of $\Delta^\calA$.

The tautological map $\pi_{\calA}$~\eqref{Eq:Taut_proj} appears in statistics.
Given (normalized) data $q \in \Delta^{\calA}$, the problem of 
\Blue{{\it maximum likelihood estimation}}
is to find the toric model parameters  (a point in $\R^{d}_{>}$) that best explain the 
data  $q$.
By Lemma 4 in \cite{DR72}, the maximum likelihood estimate for
the log-linear model $X_{\calA,w}$ is the unique point $p \in X_{\calA,w}$ such that
\[
   \pi_{\calA}(p)\ =\ \pi_{\calA}(q)
  \quad \mbox{so that}\quad p\ =\ \pi^{-1}_{\calA}(\pi_{\calA}(q)).
\]
Thus inverting the tautological projection is necessary for maximum likelihood
estimation.  

Catanese, {\it et.~al.}~\cite{CHKS} defined the \Blue{{\it maximum likelihood degree}}
of the model  $X_{\calA,w}$ to be the degree of $\pi_{\calA}^{-1}$, as an
algebraic function.
Equivalently, this is the number of complex solutions to the critical equations of the 
likelihood function, which is the degree of the tautological map $\pi_{\calA}$ from
the Zariski closure $Y_{\calA,w}$ of $X_{\calA,w}$ to $\C\P^d$.
By definition, a toric patch has rational linear precision if $\pi_{\calA}^{-1}$ is 
a rational function (an algebraic function of degree 1). 

\begin{proposition}
A toric patch $X_{\calA, w}$ has rational linear precision if and only if  $X_{\calA, w}$
has maximum likelihood degree $1$.
\end{proposition}

An important family of toric models called \Blue{{\it decomposable graphical models}}
are known to have  maximum likelihood degree 1~\cite[p.~91]{Lauritzen}. 
Therefore, these models, which are not in general B\'ezier simploids, have linear
precision. 
They typically have large $d>2$ dimension.

Darroch and Ratcliff~\cite{DR72} introduced the numerical algorithm of iterative
proportional fitting, also known as generalized iterative scaling, 
for computing the inverse of the tautological projection $\pi_{\calA}^{-1}$.
Their algorithm  requires that the data $\calA$ be in a normal, homogeneous form.

Observe first that the toric patch $X_{\calA,w}$ does not change if we translate 
all elements of $\calA$ by a fixed vector $\bb$, ($\ba\mapsto \ba+\bb$),
so we may assume that $\calA$ lies in the positive orthant $\R^d_>$.
Scaling the exponent vectors in $\calA$ by a fixed positive scalar $t\in\R_>$  also does
not change $X_{\calA,w}$ as $x\mapsto x^t$ is a homeomorphism of $\R_>$ that 
extends to a homeomorphism of $\R^d_>$.
Thus we may assume that $\calA$ lies in the standard simplex $\Delta_d$ in $\R^d$:
\[
  \Delta_d\ =\ \{ x\in\R^d\mid x\geq 0\quad\mbox{and}\quad |x|\leq 1\}\,.
\]
Lastly, we lift this to the probability simplex $\Delta^+_d\in\R^{d+1}_{\geq}$,
\[ 
   \Blue{\Delta^+_d}\ :=\ 
     \{ y\in \R^{1+d}\mid y_i\geq 0\quad\mbox{and}\quad |y|=1\}\,,
\]
by 
\[
   \calA\ni \ba \ \longmapsto\ \ba^+:=(1{-}|\ba|,\,\ba) \in \Delta^+_d\,.
\]
Since for $t\in\R_>$ and $x\in\R^d_>$,
\[
   (t,tx)^{\ba^+}\ =\ t^{1-|\ba|}(tx)^\ba\ =\ 
    t x^\ba\,,
\]
we see that replacing $\calA$ by its lifted version $\calA^+$ also does not change 
$X_{\calA,w}$.
We remark that this is just a way to homogenize the data for our problem.

Since these  are affine transformations, the property of 
affine invariance for patches shows that the affine map that is the composition 
of these transformations intertwines the original tautological projection with the new one.

We describe the algorithm of \Blue{{\it iterative proportional fitting}}, which is Theorem~1  
in~\cite{DR72}. 

\begin{proposition}\label{P:IPF}
 Suppose that $\calA\subset\Delta_d$ and $y\in\mbox{conv}(\calA)$.
 Then the sequence
\[
   \{ \bp^{(n)}\;\mid\; n=0,1,2\dotsc\} 
\]
 whose $\ba$-coordinates are defined by $p_\ba^{(0)}:=w_\ba$ and, for 
 $n\geq 0$, 
\[
   p_\ba^{(n+1)}\ :=\ p_\ba^{(n)}\cdot \frac{y^\ba}{(\pi_{\calA}(\bp^{(n)}))^\ba}\,,
\]
 converges to the unique point $\bp\in X_{\calA,w}$ such that 
 $\pi_{\calA}(\bp)=\pi_{\calA}(y)$.
\end{proposition}

We remark that if $\calA$ is not homogenized then to compute $\pi_{\calA}^{-1}(y)$ for 
$y\in \mbox{conv}(\calA)$, we first put
$\calA$ into homogeneous form $\calA^+$ using an affine map $\psi$, and then 
use iterative proportional fitting to compute
$\pi_{\calA^{+}}^{-1}(\psi(y))=\pi_{\calA}^{-1}(y)$.  
We also call this modification of the algorithm of Proposition~\ref{P:IPF}
iterative proportional fitting.
Thus iterative proportional fitting computes the inverse image of the tautological
projection, which, by Theorem~\ref{Th:LP} gives the unique parametrization of
$X_{\calA,w}$ having linear precision.

\begin{corollary}
 Iterative proportional fitting computes the unique parametrization of a toric patch having 
 linear precision. 
\end{corollary}

Little is known about the covergence of iterative proportional fitting.
We think that it is an interesting question to investigate this convergence, particularly
how it may relate to maximum likelihood degree.

\bibliographystyle{amsplain}
\bibliography{bibl}

\end{document}